\documentclass[11pt]{article} % use larger type; default would be 10pt
 
\usepackage[utf8]{inputenc} % set input encoding (not needed with XeLaTeX)
\usepackage[T1]{fontenc}
\usepackage[T1]{CJKutf8} 
\usepackage[english]{babel}
\usepackage{amsmath, amsthm, amssymb}
\usepackage{amssymb,amsfonts,textcomp}
\usepackage{supertabular}
\usepackage{hhline}
\usepackage{multicol}
\usepackage{array}
\usepackage{graphicx} 
\usepackage{color}
\usepackage{hyperref}
\usepackage{placeins}
\usepackage{algorithmicx}
\usepackage[ruled]{algorithm}
\usepackage{algpseudocode}
\usepackage{appendix}
\usepackage{footnote}
\usepackage{subfigure}
\usepackage{caption}
\usepackage{color, colortbl}

\theoremstyle{plain}

\theoremstyle{definition}
\newtheorem{definition}{Definition}

\newtheorem{open problem}{Open Problem}

\theoremstyle{remark}

\definecolor{Gray}{gray}{0.9}
\definecolor{LightCyan}{rgb}{0.88,1,1}
\definecolor{Red}{rgb}{1,0,0}
\definecolor{Orange}{rgb}{1,0.5,0}

\usepackage{geometry} % to change the page dimensions
\usepackage{fancyhdr} % This should be set AFTER setting up the page geometry
\usepackage{sectsty}
\allsectionsfont{\sffamily\mdseries\upshape} % (See the fntguide.pdf for font help)
\def\@fnsymbol#1{\ensuremath{\ifcase#1\or *\or \dagger\or \ddagger\or
   \mathsection\or \mathparagraph\or \|\or **\or \dagger\dagger
   \or \ddagger\ddagger \else\@ctrerr\fi}}

\title{Sub-problems of the $(3, 14)$ cage problem and their computer analysis}

\author{Vivek S. Nittoor\\
\small\tt vivek@nittoor.com\\
\small\tt Independent Consultant \& Researcher \footnote{formerly with the University of Tokyo}}
\date{} % Activate to display a given date or no date (if empty),

\begin{document}
\maketitle

\begin{abstract}
A $(k, g)$ graph is a graph with regular degree $k$ and girth $g$. The cage problem refers to finding the smallest $(k, g)$ graph. The $(3, 14)$ cage problem is known to be unresolved. In 2002, Exoo found a $(3, 14)$ record graph with order 384. The trivalent cage problem is restricted in this paper to the Hamiltonian bipartite class of trivalent graphs. A parameter called \textit{symmetry factor} for representing rotational symmetry is introduced in this paper. The general problem of finding a $(3, g)$ Hamiltonian bipartite graph of minimum order is further decomposed into a set of sub-problems for finding $(3, g)$ Hamiltonian bipartite graphs of minimum order for various symmetry factors. The minimum order for $(3, g)$ Hamiltonian bipartite graphs for various symmetry factors has been found using computer search. This information about sub-problems also yields useful information about non-existence of $(3, 14)$ Hamiltonian bipartite graphs between the $(3, 14)$ lower bound, 258 and the $(3, 14)$ upper bound, 384. This non-existence information partially supports the likelihood of the current $(3, 14)$ record graph indeed being the $(3, 14)$ cage.

\end{abstract}
%(Some part of the search space has been found to ....)
%(This approach could be potentially useful to completely the search the space of graphs between the lower bound, and conclusively solve the $(3, 14)$ cage problem).
%The $(3, 14)$ cage problem has not been resolved ..... lower bound .... upper bound.
%Choice of HB has been discussed in \cite{OverallPaper}.


\section{Introduction}
In $1947$, Tutte \cite{42} posed the cage problem as a problem in extremal graph theory. A $(k, g)$ graph is a graph with regular degree $k$ and girth $g$. The $(k, g)$ cage problem deals with finding a $(k, g)$ graph with minimum order. One of the well known survey papers for the cage problem is Exoo et al. 2011 \cite{Jajcaysurvey}.\\
The present paper is related to the catalog of $(3, g)$ Hamiltonian bipartite graphs (HBGs) in \cite{CatalogPaper} obtained by computer search, and its properties have been discussed in \cite{OverallPaper} and \cite{OverallPaper2}. The choice of a class of $(3, g)$ graphs for this catalog is Hamiltonian bipartite class, and this is explained in detail in \cite{OverallPaper}.\\
The $(3, 14)$ cage problem is currently unresolved. The lower bound for $(3, 14)$ graphs was improved from $256$ to $258$ by McKay et al. \cite{105} in $1998$. Quoting from \cite{OverallPaper}, ``The difficulty in finding the smallest $(3, g)$ graph in general is illustrated in the following historical example. The $(3, 14)$ vertex-transitive graph with order $406$ was found by Hoare \cite{Hoare1983} in 1981, and it was only in 2002 that a smaller $(3, 14)$ graph with order $384$ was found by Exoo \cite{Exoo} outside the vertex-transitive class.''\\  
The present paper does not improve the lower bounds for the $(3, 14)$ problem but does provide some useful information about the non-existence of $(3, 14)$ graphs. 
%for various conditions listed in Table \ref{table_non_existence4}. 
The information about non-existence of $(3, 14)$ graphs is of interest to the cage problem. The significance of this paper is that this reduces the search space of graphs of orders between the lower bound, 258 and upper bound, 384 for the $(3, 14)$ cage. \\
Symmetry factor is a parameter that has been introduced to represent rotational symmetry in \cite{OverallPaper}. The definition of symmetry factor from \cite{OverallPaper} is provided in Definition \ref{def_sym_fac_gen}.
\begin{definition}  \textbf{Symmetry factor for Hamiltonian trivalent bipartite graph}  \cite{OverallPaper}\\
\label{def_sym_fac_gen}
\label{def_sym_fac}
A Hamiltonian trivalent bipartite graph with order $2m$ is said to have symmetry factor $b \in \mathbb{N}$ if the following conditions are satisfied.
\begin{enumerate}
\item $b$ divides $m$. 
\item There exists a labelling of the vertices of the Hamiltonian trivalent graph with order $2m$ and labels $1, 2, \ldots, 2m$, such that $1 \to 2 \to \ldots \to 2m \to 1$ is a Hamiltonian cycle that satisfy the following properties.
\begin{itemize}
\item The edges that are not part of the above Hamiltonian cycle are connected as follows.
Vertex $i$ is connected to vertex $u_{i}$ for $1 \le i \le 2m$.
\item If $j \equiv i \bmod 2b$ for $1 \le j \le 2b$ and $1 \le i \le 2m$ then the following is true, $u_{i} - i  \equiv u_{j} - j \bmod 2m$.
\end{itemize}
\end{enumerate}
\end{definition}

%Cayley and vertex-transitive .sub-problems ....

%\section{First section}

%^Symmetry factor is a parameter for representing rotational symmetry in a Hamiltonian trivalent bipartite graph that has been defined in \cite{OverallPaper}. Quoting from \cite{OverallPaper}, ``Symmetry factor allows the decomposition of the problem of listing $(3, g)$ HBGs to sub-problems of listing $(3, g)$ HBGs for a range of symmetry factors.'' 

The following observations about symmetry factor are obvious from Definition \ref{def_sym_fac}.
\begin{enumerate}
\item If a HBG with order $2m$ has symmetry factor $b$, then it also has symmetry factor $m$. This case of a HBG of order $2m$ and symmetry factor $m$, is referred to as \textit{full symmetry factor}.
\item If a HBG with order $2m$ has symmetry factor $b$, then it also has symmetry factor $ab$ if $a$ is a natural number such that $ab$ divides $m$. 
\item If a HBG with order $2m$ does not exist for symmetry factor $b$, then HBG with order $2m$ does not exist for symmetry factor $a$, where $a$ is a natural number such that $a$ divides $b$.
\end{enumerate}

\FloatBarrier
%I introduce the D3 chord index notation for representing HBGs in \cite{OverallPaper}, but do not use this notation in this paper.\\
\section{$(3, 14)$ Sub-problems}
Sub-problems of the cage problems for classes of trivalent graphs such as Cayley graphs and vertex-transitive graphs have been considered by many researchers as summarized in Exoo et al. 2011 \cite{Jajcaysurvey}.\\
In the present paper, a new kind of a sub-problem is introduced, i.e., finding $(3, g)$ HBGs for various symmetry factors. These sub-problems are not disjoint.  Quoting from \cite{OverallPaper}, ``Symmetry factor allows decomposing the problem of listing $(3, g)$ HBGs for even girth $g$ for a range of orders into sub-problems of listing $(3, g)$ HBGs for a specified symmetry factor $b$ for a range of orders, and hence allows listing of $(3, g)$ HBGs for more orders.''\\
The upper bounds for $(3, 14)$ HBGs for various symmetry factors obtained from the catalog of $(3, g)$ HBGs listed in \cite{CatalogPaper} has been presented in Table \ref{table_comp314}. The lower bound for $(3, 14)$ is known to be 258 as per \cite{Jajcaysurvey}. The lower bound for a $(3, 14)$ HBG for a particular symmetry factor $b$, would be the smallest positive integer greater than or equal to 258 that is also divisible by $2b$, which we denote as $lb(3, 14, b)$.\\
There are no $(3, 14)$ HBGs of symmetry factor $3$ between 258 and 384, as shown in Table \ref{table_comp314}. For symmetry factors 4, 5, 6, the lower bound and the upper bound are equal, and hence the $(3, 14)$ sub-problems for these symmetry factors are resolved. For symmetry factors 7, 8, 9 the lower bound has been improved significantly beyond the known $(3, 14)$ lower bound, 258. Exoo's $(3, 14)$ record graph has been found on the catalog of $(3, g)$ graphs and is found for symmetry factor 8. A current lower bound for $(3, 14)$ HBGs with symmetry factor 8 is 304 as shown in Table \ref{table_comp314}. There is scope for improvement of lower and upper bound for symmetry factors $10$ to $16$.

%for sub-problems for various symmetry factor is given . The graph

\begin{table}
\centering
\caption{$(3, 14)$ sub-problems lower and upper bounds}
\label{table_comp314}
\begin{tabular}{cllllll}
\hline
Symmetry & $lb(3, 14, b)$  & Lower bound $(3, 14)$ HBG & Upper bound $(3, 14)$ HBG\\
factor $b$ &  &  for symmetry factor factor $b$ & for symmetry factor $b$ \\
\noalign{\smallskip}
\hline
\noalign{\smallskip}

\cellcolor{Red}3 &\cellcolor{Red}258 & \cellcolor{Red}900  & \cellcolor{Red}\\
\cellcolor{Gray}4 &\cellcolor{Gray}264 &\cellcolor{Gray}440 & \cellcolor{Gray}440 Figure \ref{graph_440_sym4_g14}\\
\cellcolor{Gray}5 &\cellcolor{Gray}260& \cellcolor{Gray}460 & \cellcolor{Gray}460 Figure \ref{graph_460_sym5_g14}\\
\cellcolor{Gray}6  &\cellcolor{Gray}264 & \cellcolor{Gray}456 & \cellcolor{Gray}456 Figure \ref{graph_456_sym6_g14}\\
\cellcolor{Orange}7 &\cellcolor{Orange}266 &\cellcolor{Orange}364 & \cellcolor{Orange}406 Figure \ref{graph_406_sym7_g14}\\
\cellcolor{Orange}8  &\cellcolor{Orange}272 &\cellcolor{Orange}304& \cellcolor{Orange}\textcolor{blue}{384} Figure \ref{graph_384_sym8_g14}\\
\cellcolor{Orange}9 &\cellcolor{Orange}270 &\cellcolor{Orange}288  & \cellcolor{Orange}504 Figure \ref{graph_504_sym9_g14}\\
\cellcolor{LightCyan}10 &\cellcolor{LightCyan}260 &\cellcolor{LightCyan}260 & \cellcolor{LightCyan}460  Figure \ref{graph_460_sym5_g14}\\
\cellcolor{LightCyan}11 &\cellcolor{LightCyan}264 & \cellcolor{LightCyan}264 & \cellcolor{LightCyan}506 Figure \ref{graph_506_sym11_g14}\\
\cellcolor{LightCyan}12 &\cellcolor{LightCyan}264& \cellcolor{LightCyan}264 & \cellcolor{LightCyan}456 Figure \ref{graph_456_sym6_g14}\\
\cellcolor{LightCyan}13 &\cellcolor{LightCyan}260 & \cellcolor{LightCyan}260 & \cellcolor{LightCyan}572 Figure \ref{graph_572_sym13_g14}\\
\cellcolor{LightCyan}14 &\cellcolor{LightCyan}280 & \cellcolor{LightCyan}280 & \cellcolor{LightCyan}588 Figure \ref{graph_588_sym14_g14}\\
\cellcolor{LightCyan}15 &\cellcolor{LightCyan}270 & \cellcolor{LightCyan}270 & \cellcolor{LightCyan}510 Figure \ref{graph_510_sym5_g14}\\
\cellcolor{LightCyan}16 &\cellcolor{LightCyan}288 & \cellcolor{LightCyan}288 & \cellcolor{LightCyan}\textcolor{blue}{384} Figure \ref{graph_384_sym8_g14}\\
\hline
\end{tabular}
\begin{tabular}{lclllll}
\hline\noalign{\smallskip}
Color & Significance  \\  
\noalign{\smallskip}
\hline
\noalign{\smallskip} 
\label{table_sub_problems_colors} 
 \cellcolor{Red} & Found to not exist \\
\cellcolor{Gray} & Lower bound equals upper bound \\
\cellcolor{Orange} & Lower bound improved over $lb$\\
\cellcolor{LightCyan} & Scope for potentially improving bounds  \\
\hline
\end{tabular}\\
\textcolor{blue}{384} refers to the $(3, 14)$ Record graph \\
\end{table}

\FloatBarrier

\begin{figure}[htpb]
\centering
\begin{tabular}{cccc}

%\subfigure{\includegraphics[scale=0.35]{../Cage/G14_Sym4/440_sym4.pdf}   \label{graph_440_sym4_g14}}   & 
%\subfigure{\includegraphics[scale=0.35]{../Cage/G14_Sym5_Graphs/460_Sym5.pdf}   \label{graph_460_sym5_g14}}   & 
%\subfigure{\includegraphics[scale=0.35]{../Cage/G14_Sym6_Graphs/456_Sym6.pdf}   \label{graph_456_sym6_g14}}   & 
%\subfigure{\includegraphics[scale=0.35]{../Cage/G14_Sym7_Graphs/406_Sym7.pdf}   \label{graph_406_sym7_g14}}   

\subfigure{\includegraphics[scale=0.35]{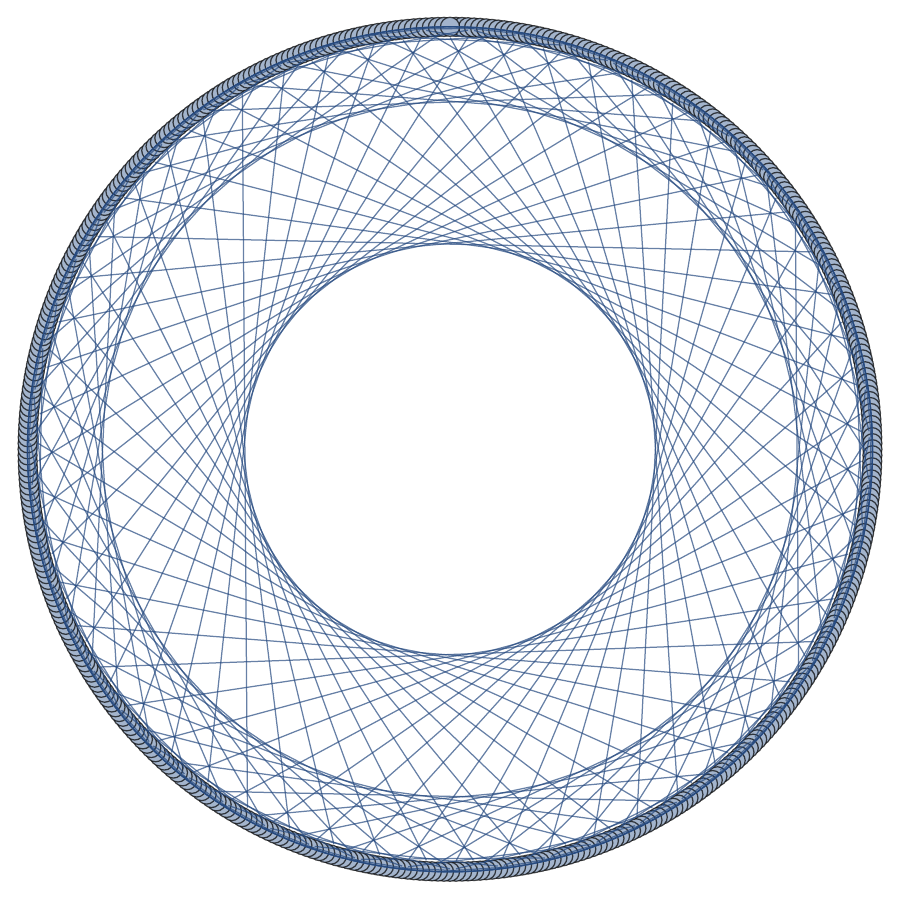}   \label{graph_440_sym4_g14}}   & 
\subfigure{\includegraphics[scale=0.35]{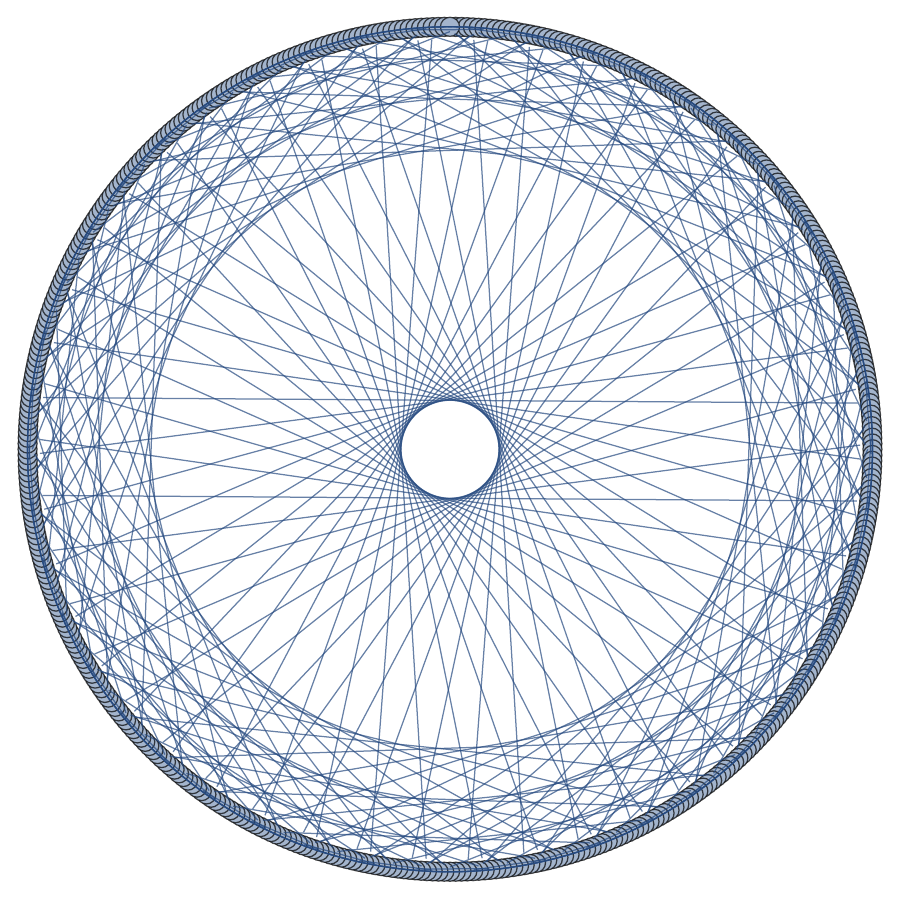}   \label{graph_460_sym5_g14}}   & 
\subfigure{\includegraphics[scale=0.35]{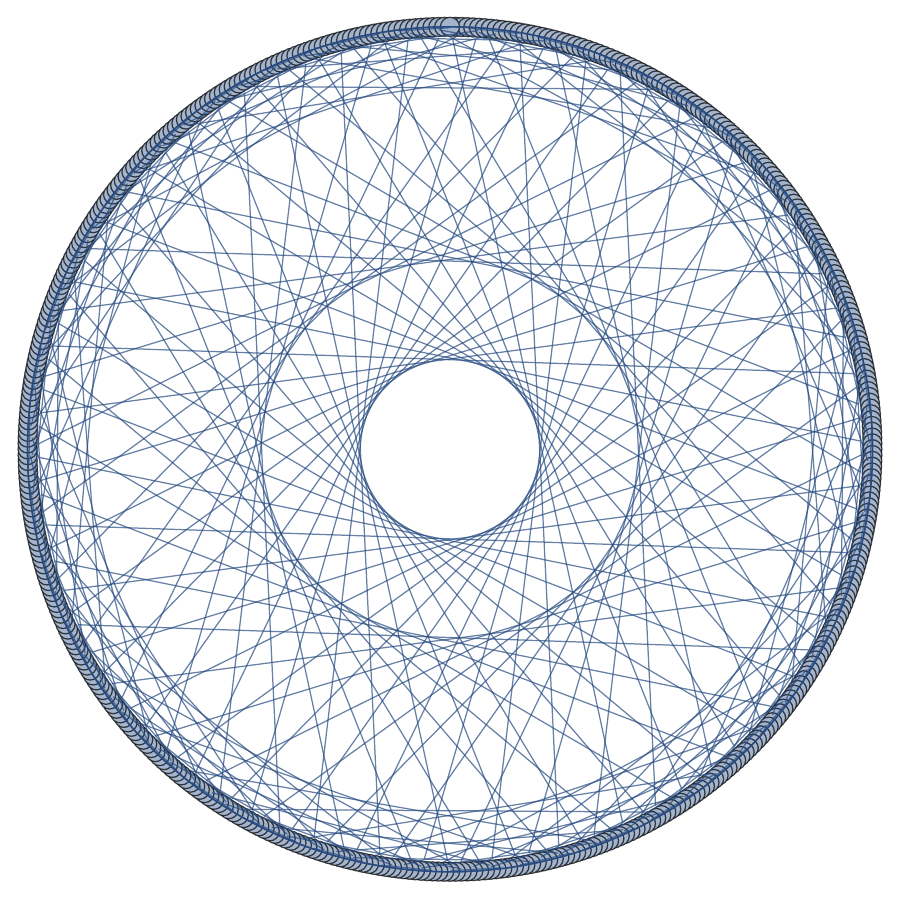}   \label{graph_456_sym6_g14}}   & 
\subfigure{\includegraphics[scale=0.35]{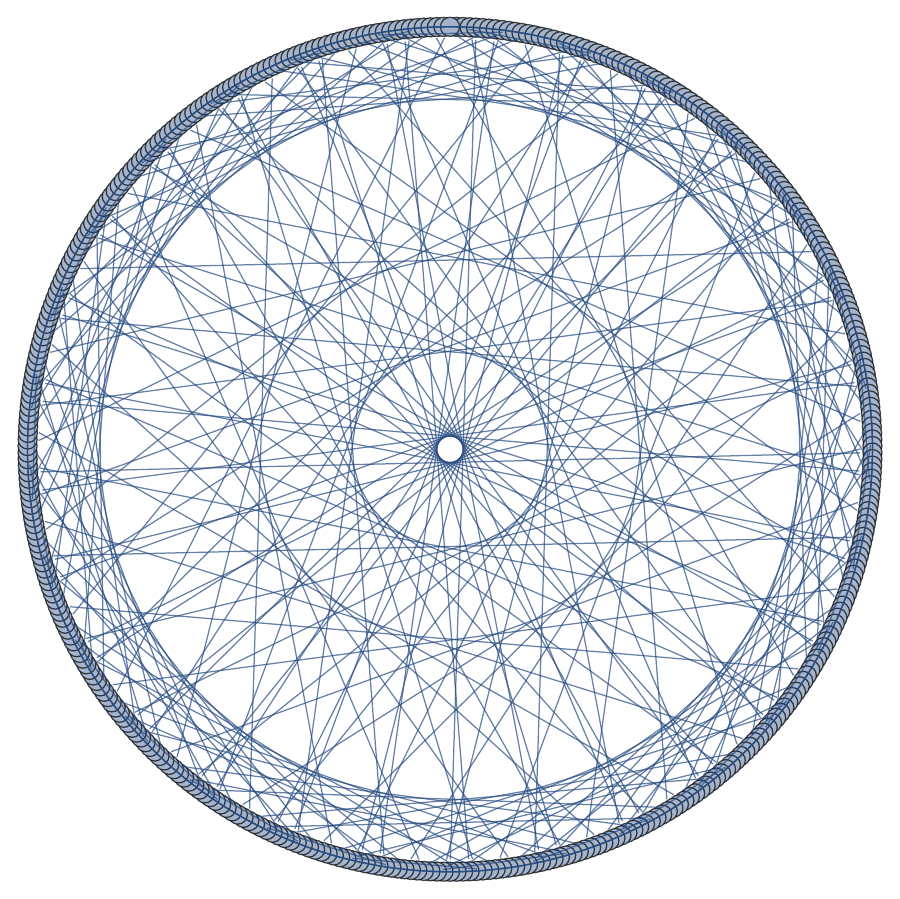}   \label{graph_406_sym7_g14}}   

\\
(a) $440$ sym $4, 20, 44$ & (b) $460$ sym $5, 10$ & (c) $456$ sym $6, 12$ & (d) $406$ sym $7$
\\  
%\subfigure{\includegraphics[scale=0.35]{../Cage/384.pdf}   \label{graph_384_sym8_g14}}  &
%\subfigure{\includegraphics[scale=0.35]{../Cage/G14_Sym9_Graphs/504_Sym9.pdf}   \label{graph_504_sym9_g14}} &
%\subfigure{\includegraphics[scale=0.35]{../Cage/G14_Sym11_Graphs/506_Sym11.pdf}   \label{graph_506_sym11_g14}} &
%\subfigure{\includegraphics[scale=0.35]{../Cage/G14_Sym13_Graphs/572_Sym13.pdf}   \label{graph_572_sym13_g14}}

\subfigure{\includegraphics[scale=0.35]{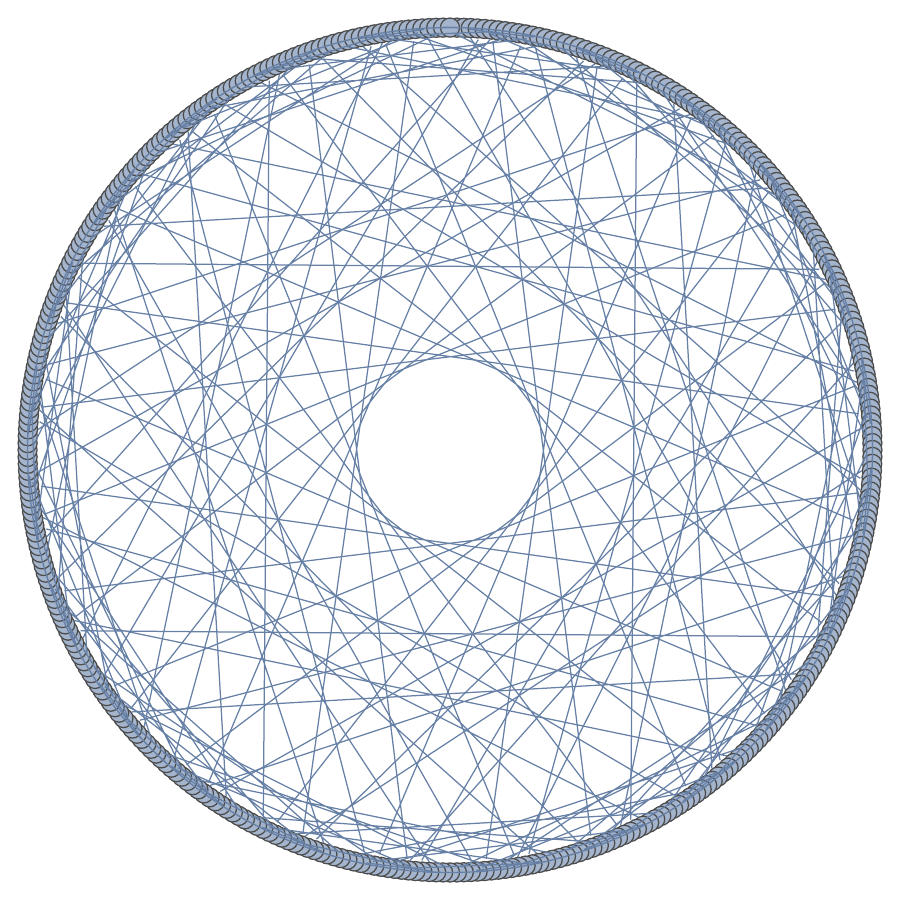}   \label{graph_384_sym8_g14}}  &
\subfigure{\includegraphics[scale=0.35]{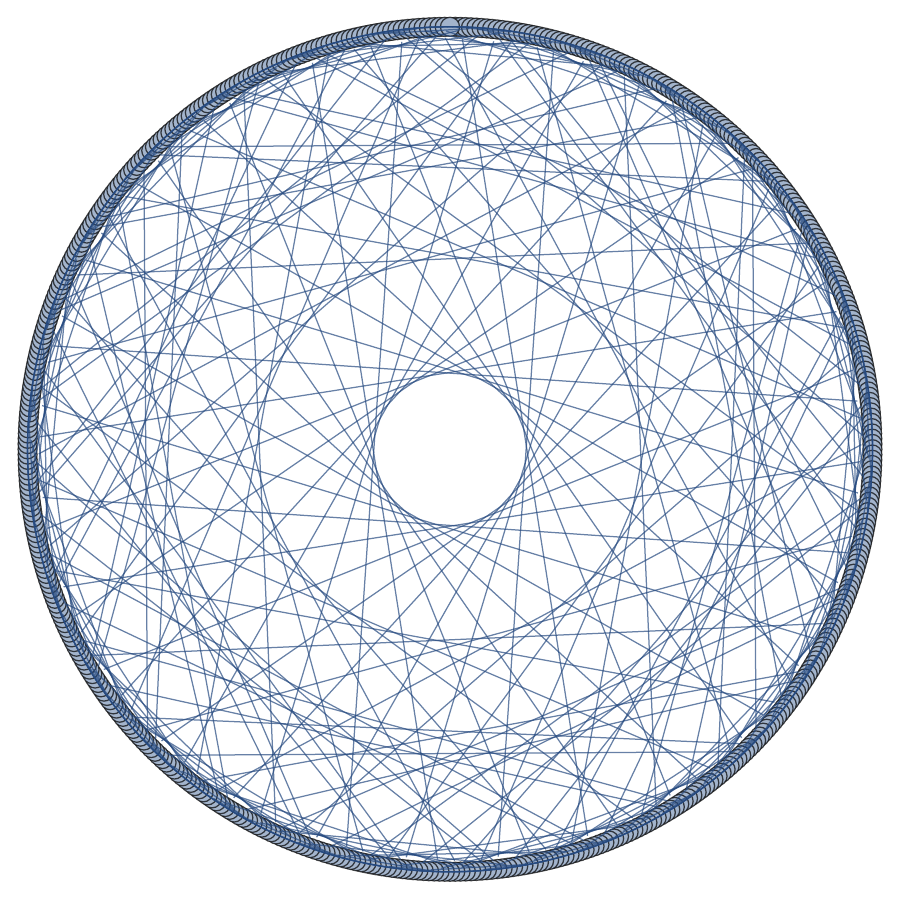}   \label{graph_504_sym9_g14}} &
\subfigure{\includegraphics[scale=0.35]{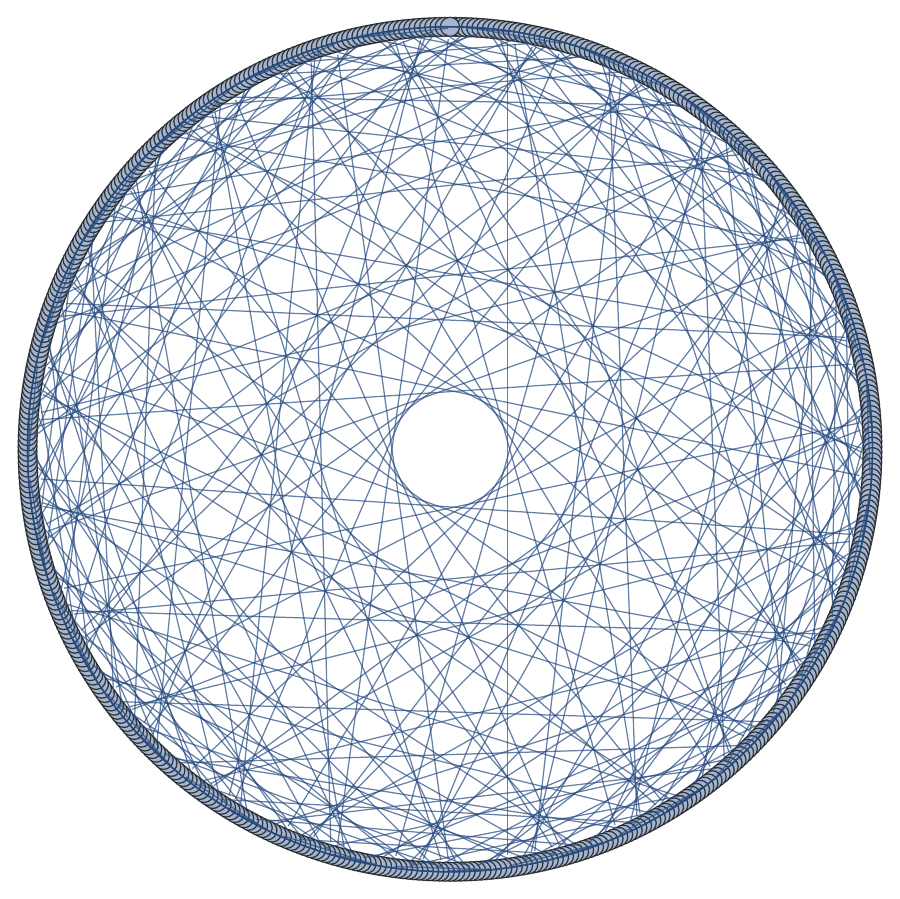}   \label{graph_506_sym11_g14}} &
\subfigure{\includegraphics[scale=0.35]{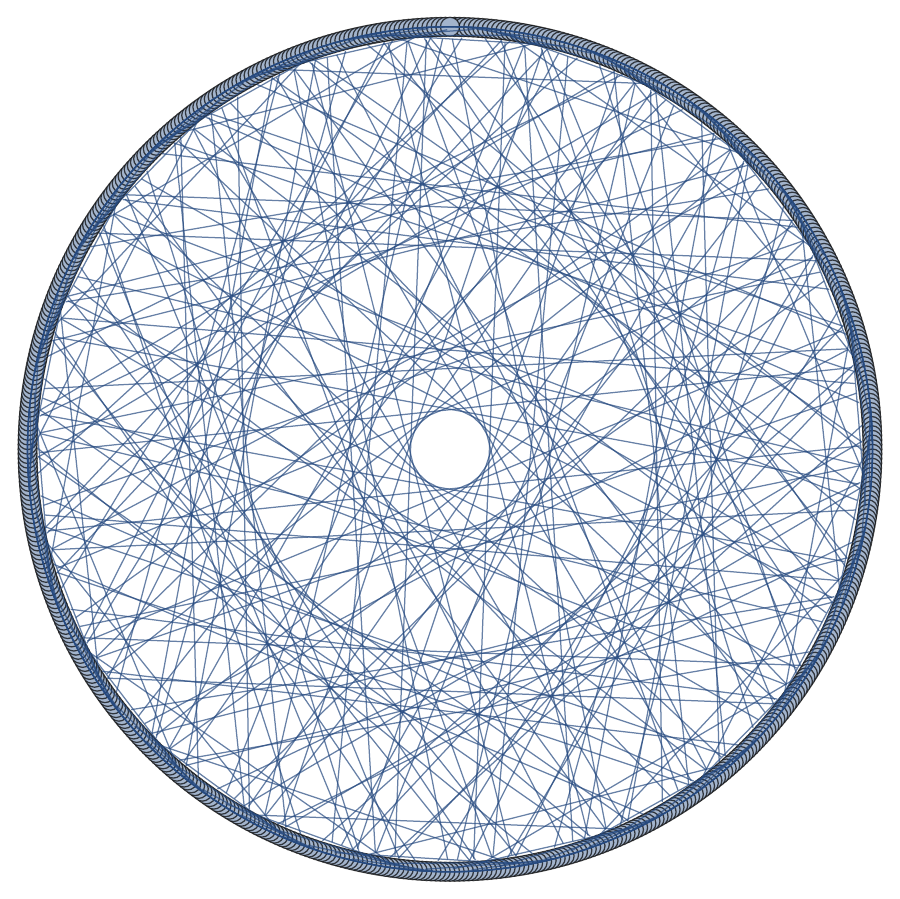}   \label{graph_572_sym13_g14}}

\\

(e) $384$ sym $8, 16, 24, 48, 96$ & (f) $504$ sym $9, 18, 36$ & (g) $506$ sym $11$  &  (h) $572$ sym $13, 26, 52, 143$ \\

%\subfigure{\includegraphics[scale=0.35]{../Cage/G14_Sym14_Graphs/588_Sym14.pdf}   \label{graph_588_sym14_g14}} &
%\subfigure{\includegraphics[scale=0.35]{../Cage/G14_Sym5_Graphs/510_Sym5.pdf}   \label{graph_510_sym5_g14}}

\subfigure{\includegraphics[scale=0.35]{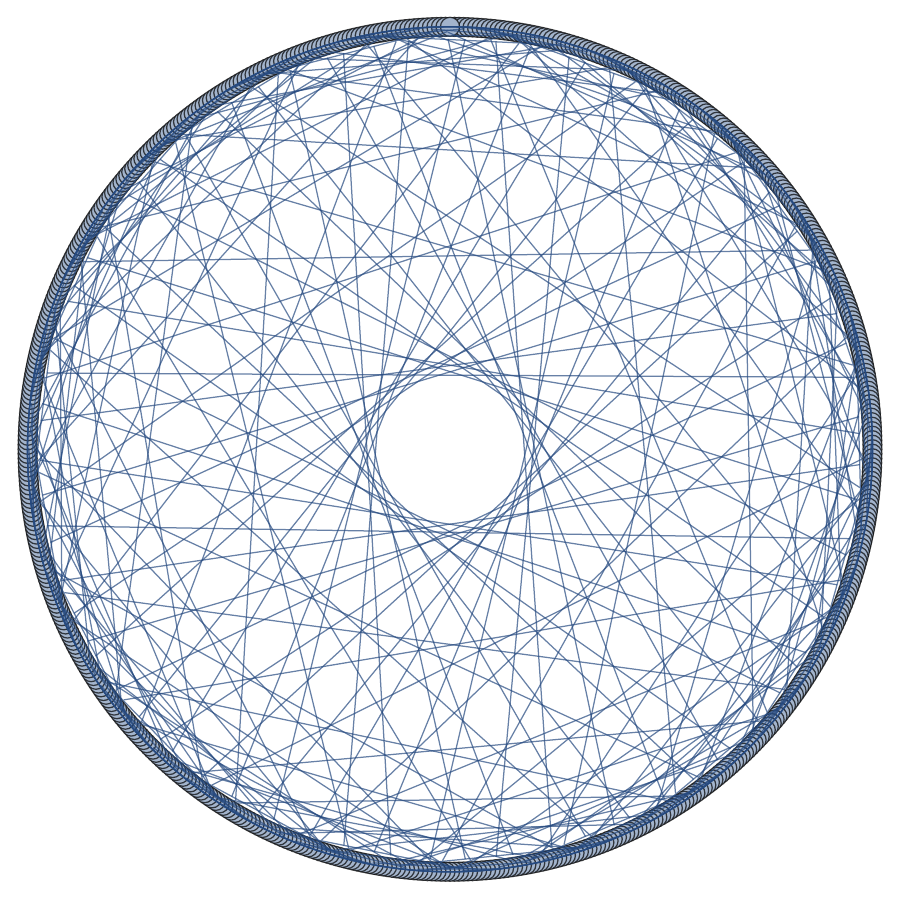}   \label{graph_588_sym14_g14}} &
\subfigure{\includegraphics[scale=0.35]{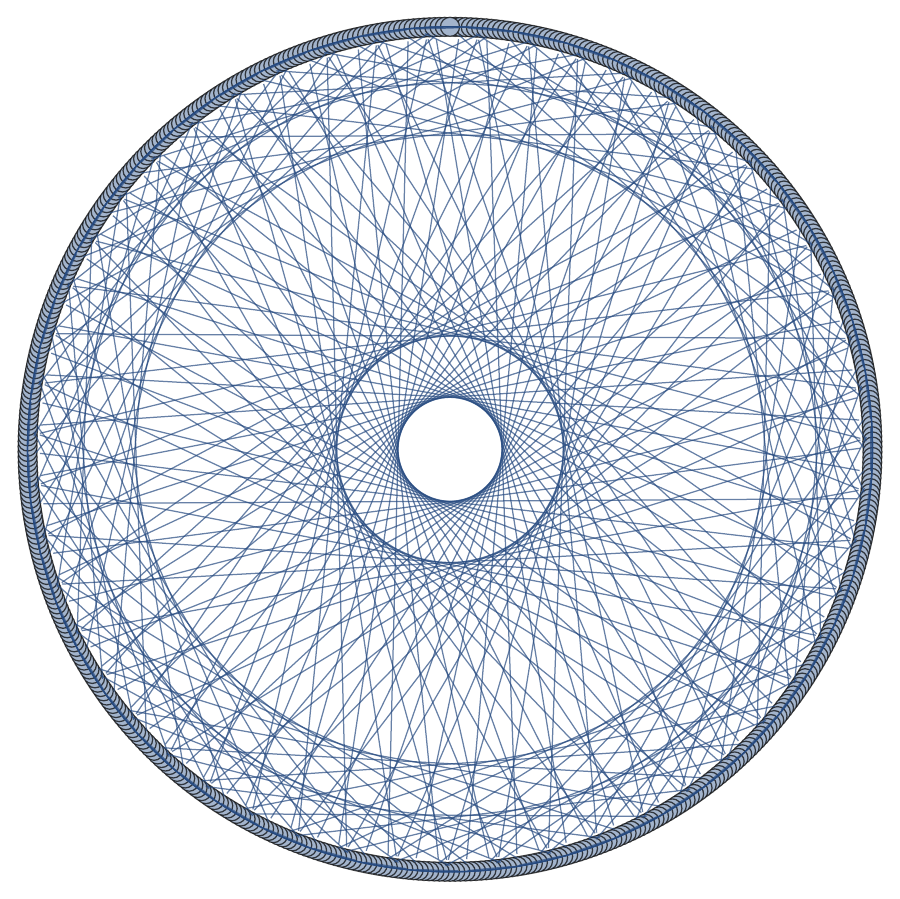}   \label{graph_510_sym5_g14}}

\\

(i) $588$ sym $14, 28, 49, 98$ & (j) $510$ sym $15$, (5) \\

\end{tabular}
\caption{Upper bounds for $(3, 14)$ HBGs for various symmetry factors from Catalog \cite{CatalogPaper}}
\label{table_3_14_Min_graphs}
\end{figure}

In Figure \ref{table_3_14_Min_graphs}, $(3, 14)$ HBGs that are upper bounds for various symmetry factors have been shown, using the notation of order of the graph followed by ``sym'' and the symmetry factors for the HBG. \textit{full symmetry factor} is not mentioned in each these cases since it is obvious. \\ 
In cases of symmetry factors where the upper bounds are equal to the lower bounds in Table \ref{table_comp314},  the problem of finding the minimum $(3, g)$ HBG with that particular symmetry factor is considered to be resolved, and the upper bound for that particular symmetry factor is refered to the minimum $(3, g)$ HBG for that particular symmetry factor.
The smallest $(3, 14)$ HBG with symmetry factors $4, 20, 44$ is of order $440$ and is shown in Figure \ref{graph_440_sym4_g14} as ``440 sym 4, 20, 44''. The smallest $(3, 14)$ HBG with symmetry factors $5$ and $10$ is of order $460$ and is shown in Figure \ref{graph_460_sym5_g14}. The smallest $(3, 14)$ HBG with symmetry factors $6$ and $12$ is of order $456$ and is shown in Figure \ref{graph_456_sym6_g14}. The smallest $(3, 14)$ HBG with symmetry factor $7$ is of order $406$ and is shown in Figure \ref{graph_406_sym7_g14}.\\
The smallest known $(3, 14)$ HBG with symmetry factors $8, 16, 24, 48, 96$ is of order $384$ and is shown in Figure \ref{graph_384_sym8_g14}, and is isomorphic to the $(3, 14)$ record graph found by Exoo \cite{Exoo} 2002.\\
The smallest known $(3, 14)$ HBG with symmetry factors $9, 18, 36$ is of order $504$ and is shown in Figure \ref{graph_504_sym9_g14}. The smallest known $(3, 14)$ HBG with symmetry factors $11$ is of order $506$ and is shown in Figure \ref{graph_506_sym11_g14}. The smallest known $(3, 14)$ HBG with symmetry factors $13, 26, 52, 143$ is of order $572$ and is shown in Figure \ref{graph_572_sym13_g14}. The smallest known $(3, 14)$ HBG with symmetry factors $15$ is of order $510$ and is shown in Figure \ref{graph_510_sym5_g14}. This HBG also has symmetry factor $5$ but is not the smallest $(3, 14)$ HBG of symmetry factor $5$, and hence 5 is shown in brackets in Figure \ref{graph_510_sym5_g14}. The smallest known $(3, 14)$ HBG with symmetry factors $14, 28, 49, 98$ is of order $588$ and is shown in Figure \ref{graph_588_sym14_g14}. 

%\FloatBarrier

There is currently no other known source of information about non-existence of $(3, 14)$ graphs between orders 258 and 384. Non-existence of $(3, 14)$ HBGs for several orders and symmetry factors is implied from Table \ref{table_comp314}. For the convenience of the reader, this information about non-existence  for $(3, 14)$ HBGs for various symmetry factors has been summarised in Table \ref{table_non_existence4}. More detailed non-existence information has been provided in \cite{OverallPaper2}. The additional information in Table \ref{table_non_existence4} compared to Table  \ref{table_comp314} is the non-existence of $(3, 14)$ HBGs with symmetry factor 19 and order 380.\\
One interesting observation is that the smallest $(3, 14)$ HBG with symmetry factor 4 is 440, but a $(3, 14)$ HBG with symmetry factor 4 and order 456 does not exist as shown in Table 
\ref{table_non_existence4}. This phenomenon is referred to as \textit{non-monotonicity} and this is discussed in \cite{OverallPaper}.\\

\begin{table}
\centering
\caption{Non-existence of $(3, 14)$ HBGs for various symmetry factors for the following orders}
\label{table_non_existence4} 
\begin{tabular}{clclllll}
\hline\noalign{\smallskip}
%Symmetry & $(3, 14)$ Non-existence \\  
Symmetry factor $b$ & Orders for non-existence of $(3, 14)$ HBG\\
 & with symmetry factor $b$ \\
\noalign{\smallskip}
\hline
\noalign{\smallskip} 
$3$ & 258, 264, 270, 276, 282, 288, 294, 300, 306, 312, 318, 324, 330, 336, 342, \\
& 348, 354, 360, 366, 372, 378, 384, ....	\\
$4$ & 264, 272, 280, 288, 296, 304, 312, 320, 328, 336, 344, 352, 360, 368, 376, 384,   \\  
 & 392, 400, 408, 416, 424, 432, 456 \\  
$5$ & 260, 270, 280, 290, 300, 310, 320, 330, 340, 350, 360, 370, 380, 390, 400, \\   
 & 410, 420, 430, 440, 450, 470, 480 \\   
$6$ & 264, 276, 288, 300, 312, 324, 336, 348, 360, 372, 384\\ & 396, 408, 420, 432, 444  \\  
$7$ &266, 280, 294, 308, 322, 336, 350 \\  
$8$ & 272, 288 \\  
$9$ &  270 \\
$19$ & 380 \\
\hline
\end{tabular}
\end{table}

\FloatBarrier

\section{Conclusion}
The present paper hence reduces the search space of possible $(3, 14)$ graphs between orders 258 and 384, for the possible existence of a $(3, 14)$ graph with order within the range, and hence partially supports the likelihood of Exoo's $(3, 14)$ graph with order $384$ being the $(3, 14)$ cage.

%\section*{Acknowledgement}
%\textbf{Acknowledgments}\\ 
%This research did not receive any specific grant from funding agencies in the public, commercial, or not-for-profit sectors.

%-------------------
\bibliographystyle{plain} 
\bibliography{vivek}
%%%%%%%%%%%%%%%%      

\end{document}